\documentclass{amsart}
\usepackage{amsmath, graphicx}

\newcommand\cC{\mathcal{C}}
\newcommand\cL{\mathcal{L}}
\newcommand\cO{\mathcal{O}}

\newcommand\cX{\mathcal{X}}

\renewcommand\d{\mathrm{d}}
\newcommand\e{\mathrm{e}}
\renewcommand\i{\mathrm{i}}

\newcommand\C{\mathbf{C}}
\newcommand\R{\mathbf{R}}

\newcommand\oA{\bar{A}}
\newcommand\oE{\bar{E}}
\newcommand\oL{\bar{L}}
\newcommand\ou{\bar{u}}
\newcommand\ov{\bar{v}}
\newcommand\oy{\bar{y}}
\newcommand\oOm{\bar{\Omega}}

\newcommand\hu{\hat{u}}

\renewcommand\Re{\operatorname{Re}}
\renewcommand\Im{\operatorname{Im}}

\let\al\alpha
\let\be\beta
\let\ga\gamma

\let\eps\varepsilon
\let\ka\kappa
\let\la\lambda
\let\La\Lambda
\let\si\sigma
\let\phi\varphi
\let\Om\Omega

\newcommand\sipt{\sigma_{\text{pt}}}
\newcommand\siess{\sigma_{\text{ess}}}

\begin{document}

\title[Evaluating the Evans function]%
      {Evaluating the Evans function:\\Order reduction in numerical methods}
\author{Simon Malham}
\address{Mathematics Department\\Heriot-Watt University\\Edinburgh
  EH14 4AS\\United Kingdom}
\email{simonm@ma.hw.ac.uk}
\author{Jitse Niesen}
\address{Mathematics Department\\Heriot-Watt University\\Edinburgh
  EH14 4AS\\United Kingdom. (Current address: Mathematics Department
  \\ La Trobe University \\ Victoria 3086 \\ Australia)}
\email{j.niesen@latrobe.edu.au}
\thanks{This work was supported by EPSRC First Grant GR/S22134/01.}

\begin{abstract}
  We consider the numerical evaluation of the Evans function, a
  Wronskian-like determinant that arises in the study of the stability
  of travelling waves. Constructing the Evans function involves
  matching the solutions of a linear ordinary differential equation
  depending on the spectral parameter. The problem becomes stiff as
  the spectral parameter grows. Consequently, the Gauss--Legendre
  method has previously been used for such problems; however more
  recently, methods based on the Magnus expansion have been
  proposed. Here we extensively examine the stiff regime for a general
  scalar Schr\"odinger operator. We show that although the
  fourth-order Magnus method suffers from order reduction, a
  fortunate cancellation when computing the Evans matching function
  means that  fourth-order convergence in the end result is preserved.
  The Gauss--Legendre method does not suffer from order
  reduction, but it does not experience the cancellation either, and
  thus it has the same order of convergence in the end result. Finally
  we discuss the relative merits of both methods as spectral tools.
\end{abstract}

\maketitle

\section{Introduction}
\label{sec:intro}

Many partial differential equations admit travelling wave solutions;
these are solutions that move at a constant speed without changing
their shape. Such travelling waves occur in many fields, including
biology, chemistry, fluid dynamics, and optics. It is often important
to know whether a given travelling wave is stable: does it persist
under small perturbations? A major step towards determining the
stability of a travelling wave is to locate the spectrum of the
linearization of the differential operator about the travelling
wave. Evans~\cite{evans:nerve*3} considers a shooting and matching
method for this task. Evans introduced a function to measure the
mismatch for a specific class of reaction--diffusion equations. This
function was called the \emph{Evans function} by Alexander, Gardner
and Jones~\cite{alexander.gardner.ea:topological}, who generalized its
definition considerably. Since then, the Evans function has been used
frequently for stability analysis; see, for example,
\cite{afendikov.bridges:instability, aparicio.malham.ea:numerical,
bridges.derks.ea:stability, brin:numerical,
gubernov.mercer.ea:evans*1} for numerical computations employing the
Evans function and \cite{facao.parker:stability,
kapitula.sandstede:edge, pego.weinstein:eigenvalues,
swinton.elgin:stability, terman:stability} for an analytic
treatment. The review paper by Sandstede~\cite{sandstede:stability}
gives an excellent overview of the field.

The Evans function is a function of one argument, the spectral
parameter~$\la$, and zeros of the Evans function correspond to
eigenvalues of the corresponding operator. Hence, one can get
information about the spectrum by finding zeros of the Evans function,
either analytically or numerically. Our focus here is on the numerical
approach.

The main part in the numerical evaluation of the Evans function is the
solution of a linear ordinary differential equation depending on the
spectral parameter~$\la$.  This is often done with an off-the-shelf
integrator using an explicit Runge--Kutta method. Afendikov and
Bridges~\cite{afendikov.bridges:instability} noticed that when $\la$
grows, the problem may become \emph{stiff}, and therefore they use a
Gauss--Legendre method.  Recently, Aparicio, Malham and
Oliver~\cite{aparicio.malham.ea:numerical} proposed a new procedure
based on the Magnus expansion, building on the work of
Moan~\cite{moan:efficient} and Greenberg and
Marletta~\cite{greenberg.marletta:numerical} who used the Magnus
expansion to solve Sturm--Liouville problems.  Aparicio \emph{et al.}\
noticed that the Magnus method suffers from \emph{order reduction} in
the stiff regime. This means that the fourth-order integrators whose
global error should scale like $h^4$ when the step size~$h$ is small,
instead converge more slowly.  They analyzed this phenomenon in a
modified Airy equation using the WKB-method. However, they restricted
themselves to those values of~$\la$ which correspond to the essential
spectrum of the linearized differential operator.

The current paper continues the analysis of the Magnus method in the
context of Evans function evaluations. We concentrate on scalar
Schr\"odinger operators to simplify the analysis. Other methods based
on a transformation to Pr\"ufer variables~\cite{pryce:numerical}
probably perform better in the scalar setting, but these methods
cannot be used unchanged in the non-self-adjoint case where our
interest lies.

We present another approach to the analysis of the Magnus method based
on a power series expansion, which is valid for values of~$\la$
outside the essential spectrum. We will show that the Magnus method
also suffers from order reduction in this regime.  Specifically, the
relative local error is of order~$\la^{-1/2}h^2$ as $h\to0$ with
$|\la|^{1/2}h \gg 1$. However, there are two subsequent important
observations. Firstly, when going from the local to the global error,
one does not lose a factor of~$h$ (as usually), but the global error
is also of order~$\la^{-1/2} h^2$. Secondly, the order reduction
disappears completely when we evaluate the matching condition: the
relative error in the Evans function is of order~$\la^{-1/2} h^4$,
thus quartic in the step size, just as one would expect from a
fourth-order method. Since useful asymptotic estimates invoke an
order~$\la^{-1}$ error, at best, our numerical schemes even with order
reduction prove a useful spectral tool in the regime $|\la| \ll
h^{-8}$.

The phenomenon of order reduction was discovered for implicit
Runge--Kutta methods by Prothero and
Robinson~\cite{prothero.robinson:on}. Nowadays, it is understood within
the framework of B-convergence (see for
instance~\cite[\S{}IV.15]{hairer.wanner:solving}). The stability of
Magnus methods has been analyzed for highly-oscillatory equation by
Iserles~\cite{iserles:on*1}, for Schr\"odinger equations by Hochbruck
and Lubich~\cite{hochbruck.lubich:on*1}, and for parabolic equations
by Gonz\'alez, Ostermann and
Thalhammer~\cite{gonzalez.ostermann.ea:second-order}. 
Unfortunately, these results cannot yet be fitted into a general
theory~\cite{iserles.munthe-kaas.ea:lie-group}. The present paper can
also be viewed as a contribution to this research.

We also present an analysis of the fourth-order Gauss--Legendre
method.  We show that the relative error committed by this method does
not contain a term of order $\la^{-1/2} h^2$, but only smaller
terms. Furthermore, the error decreases even further when evaluating
the Evans function. As explained in more detail later, this is due to
an effect similar to the one which makes the trapezoidal rule very
efficient for the quadrature of periodic functions.

The contents of this paper are as follows. In the next section, we
define the Evans function and we give an asymptotic expression for the
Evans function in the scalar case when the spectral parameter~$\la$ is
large in modulus and outside the essential spectrum.  We then define
the Magnus method in Section~\ref{sec:magnus}. We show that the Magnus
method suffers from order reduction and we compute the error when
evaluating the Evans function. We repeat the computation for the
Gauss--Legendre method in the next section. The analysis is
corraborated by numerical experiments in Section~\ref{sec:num}. In the
final section, we compare our results with those of Aparicio, Malham
and Oliver~\cite{aparicio.malham.ea:numerical} and we discuss the
stability of the Magnus method in general. More details of the
intricate calculations presented in
Sections~\ref{sec:evans}--\ref{sec:gl} can be found in the technical
report~\cite{niesen:on}.

\section{The Evans function}
\label{sec:evans}

We are interested in homogeneous reaction--diffusion equations on an
unbounded one-dimensional domain. Such equations have the form
\begin{equation}
\label{rd}
u_t = K u_{xx} + f(u),
\end{equation}
where $K$ is an $n$-by-$n$ diagonal matrix with positive entries (the
diffusion coefficients) and the unknown $u$ is a function of~$t$
and~$x$. The function $f:\R^n\to\R^n$ describes the reaction term; we
assume that $f$ is sufficiently smooth.

A travelling wave solution has the form $u(x,t) = \hu(\xi)$ with $\xi
= x-ct$ where $c$ is the wave speed---see for example Kolmogorov,
Petrovsky and Piskunov~\cite{kolmogorov.petrovsky.ea:etude}.  We
assume that a travelling wave solution for the equation is known, at
least numerically. Furthermore, we assume that $\hu$ is constant at
infinity, meaning that the limits $\hu_\pm = \lim_{\xi\to\pm\infty}
\hu(\xi)$ exist (in fact, we will need later that additionally, the
derivatives $\hu^{(p)}$ vanish at infinity for $p=1,2,\ldots$). Such a
wave is called a \emph{pulse} (if $\hu_+ = \hu_-$) or a \emph{front}
(if $\hu_+ \neq \hu_-$).

To study the stability of the travelling wave, we linearize~\eqref{rd}
about the wave and write the result in the $(\xi,t)$~coordinate system
which moves with the same speed as the travelling wave. This yields
\begin{equation}
\label{rdlin}
u_t = K u_{\xi\xi} + c u_\xi + Df(\hu)\,u.
\end{equation}
Define the operator~$\cL$ by $\cL(U) = KU'' + cU' + Df(\hu)\,U$, where
$U:\R\to\C^n$. Its spectrum determines whether~\eqref{rdlin} has
solutions of the form $u(\xi,t) = \e^{\la t} U(\xi)$. The stability of
the travelling wave~$\hu$ can be deduced from the location of the
spectrum of~$\cL$.

The spectrum of~$\cL$ can be divided in two parts: the point
spectrum~$\sipt(\cL)$, consisting of those $\la\in\si(\cL)$ for which
$\cL-\la I$ is Fredholm of index zero,%
\footnote{An operator is Fredholm with index zero if its range is
  closed and the dimension of the null space equals the codimension of
  the range.}  
and the essential spectrum~$\siess(\cL)$, which contains the rest of
the spectrum. We assume that $\siess(\cL)$ is contained in the left
half-plane $\{ z\in\C : \Re z \le 0 \}$. This means that the spectral
stability is determined by the position of the eigenvalues
$\la\in\sipt(\cL)$.

We now introduce the Evans function, which is a tool for locating
these eigenvalues. We rewrite the eigenvalue equation $\cL(U) = \la U$
as the first-order differential equation
\begin{subequations}
\label{bvp}
\begin{equation}
\label{bvp1}
\frac{\d y}{\d\xi} = A(\xi;\la)\,y,
\end{equation}
where $y:\R\to\C^{2n}$ and the matrix~$A$ is given by
\begin{equation}
\label{bvp2}
A(\xi;\la) = \begin{bmatrix}
0 & I \\ K^{-1} \bigl( \la I - Df(\hu(\xi)) \bigr) & -cK^{-1}
\end{bmatrix}.
\end{equation}
\end{subequations}
Since the spectral problem~\eqref{bvp1} is a linear equation, its
solutions form a linear space of dimension~$2n$. Define $E^-(\la)$ to
be the subspace of solutions~$y$ satisfying the boundary condition
$y(\xi) \to 0$ as $\xi \to -\infty$. Similarly, $E^+(\la)$~denotes the
subspace with $y(\xi) \to 0$ as $\xi \to +\infty$.  Any eigenfunction
must satisfy both boundary conditions and hence lie in the
intersection of $E^-(\la)$ and~$E^+(\la)$.

For all $\la\not\in\siess(\cL)$, we have 
$$
\dim E^-(\la) + \dim E^+(\la) = 2n.
$$
Choose a basis $y_1({{}\cdot{}}; \la), \ldots, y_k({{}\cdot{}}; \la)$
of~$E^-(\la)$, where $k = \dim E^-(\la)$, and a basis
$y_{k+1}({{}\cdot{}}; \la), \ldots, y_{2n}({{}\cdot{}}; \la)$
of~$E^+(\la)$. We can assemble these basis vectors, evaluated at an
arbitrary point, say $\xi=0$, in the $2n$-by-$2n$ matrix
$$
\Bigl[ \, y_1(0;\la) \,\, \ldots \,\, y_k(0;\la) \,\, y_{k+1}(0;\la)
\,\, \ldots \,\, y_{2n}(0;\la) \, \Bigr].
$$
The Evans function, denoted~$D(\la)$, is defined to be the determinant
of this matrix. If the determinant vanishes, then the $y_i$ are
linearly dependent, which implies that the spaces~$E^-(\la)$
and~$E^+(\la)$ have a nontrivial intersection, and this intersection
contains the eigenfunctions of~\eqref{bvp1}. Therefore,
$D(\la)=0$ if and only if $\la\in\sipt(\cL)$.

Let $\cC$ denote the connected component of $\C\setminus\siess(\cL)$
containing the right half-plane. We can choose the basis vectors~$y_i$
to be analytic functions of~$\la$ in the region~$\cC$. The Evans
function will then also be analytic in~$\cC$ and the order of its
zeros corresponds with the multiplicity of the eigenvalues
of~$\cL$.

More details on the Evans function and the stability of travelling
waves can be found in the landmark paper by Alexander, Gardner and
Jones~\cite{alexander.gardner.ea:topological} and the review article
by Sandstede~\cite{sandstede:stability}.

\subsection{The Evans function near infinity}

We are interested in the behaviour of~$D(\la)$ and numerical
approximations to~$D(\la)$ as $|\la| \to \infty$, because experiments
show an unexpected deterioration of the approximations in this
limit~\cite{aparicio.malham.ea:numerical}. For simplicity, we will
restrict ourselves to \emph{scalar} reaction--diffusion equations,
i.e., we assume that $n=1$. However, it is expected that the methods
of analysis presented in this paper also apply to the nonscalar case,
though the computations will obviously be more involved.

We may assume without loss of generality that the diffusion
coefficient is~1, so that the partial differential equation reads 
$$
u_t = u_{xx} + f(u).
$$ 
The corresponding eigenvalue problem~\eqref{bvp} in this case is
\begin{subequations}
\label{bvp3}
\begin{equation}
\label{bvp3a} 
\frac{\d y}{\d\xi} = A(\xi;\la)\,y, 
\end{equation}
where
\begin{equation}
\label{bvp3b}
A(\xi;\la) = \begin{bmatrix}
0 & 1 \\ \la - f'(\hu(\xi)) & -c
\end{bmatrix}.
\end{equation}
\end{subequations}
The limits of~$A$ as $\xi\to\pm\infty$ are given by 
\begin{equation}
\label{apm}
A_\pm(\la) = \begin{bmatrix}
  0 & 1 \\ \la - f'(\hu_\pm) & -c 
\end{bmatrix}.
\end{equation}
Furthermore, the eigenvalues of~$A_-(\la)$ are 
\begin{equation}
\label{mupm}
\mu_-^{[1]}, \mu_-^{[2]} = \frac12 \left( -c \pm \sqrt{c^2 + 4(\la-f'(\hu_-))}
\right).
\end{equation}
To avoid any confusion between the eigenvalues of the differential
operator~$\cL$, which form the point spectrum that we want to compute,
and the eigenvalues of the matrices~$A_\pm(\la)$, we call the latter
\emph{spatial eigenvalues}.

One of the spatial eigenvalues is purely imaginary if $\la$ lies on
the parabolic curve given by
$$
\ga_- = \bigl\{\, -s^2 + f'(\hu_-) + \i cs : s\in\R \,\bigr\}.
$$
The curve~$\ga_-$ is part of the essential spectrum.  If $\la$ lies to
the right of~$\ga_-$, then the spatial eigenvalues~$\mu_-^{[1]}$
and~$\mu_-^{[2]}$ have positive and negative real parts, respectively.

The limit $\xi\to+\infty$ is treated in the same manner and leads to
the curve~$\ga_+$. The region~$\cC$ on which the Evans function is
defined is the part of the complex plane to the right of
$\ga_-\cup\ga_+$. 

We assume henceforth that $\la\in\cC$. We compute the asymptotic
behaviour of the Evans function as $|\la|\to\infty$ using a different
approach to that outlined in Alexander, Gardner and
Jones~\cite[\S5B]{alexander.gardner.ea:topological}, extending the
approximation to further higher order corrections. We start with the
solution~$y$ of~\eqref{bvp3} satisfying $y(\xi)\to0$ as
$\xi\to-\infty$. The matrix~$A$ in~\eqref{bvp3b} goes to~$A_-$ as
defined in~\eqref{apm} in this limit, and the eigenvalues of~$A_-$ are
given in~\eqref{mupm}, with corresponding eigenvectors
$(1,\mu_-^{[1]})^\top$ and~$(1,\mu_-^{[2]})^\top$. This suggests
writing~$y$ as
\begin{subequations}
\label{oy} 
\begin{equation}
y(\xi) = \exp(\mu_-^{[1]}\xi) \left( 
  \ou(\xi) \begin{bmatrix} 1 \\ \mu_-^{[1]} \end{bmatrix}
  + \ov(\xi) \begin{bmatrix} 1 \\ \mu_-^{[2]} \end{bmatrix}
\right)  = \exp(\mu_-^{[1]}\xi) \, B \, \oy(\xi) 
\end{equation}
where
\begin{equation}
\oy = \begin{bmatrix} \ou \\ \ov \end{bmatrix} \text{ and } 
B = \begin{bmatrix} 1 & 1 \\ \mu_-^{[1]} & \mu_-^{[2]} \end{bmatrix}. 
\end{equation}
\end{subequations}
The vector $\oy$ satisfies the linear differential equation
\begin{equation}
\label{bvpo}
\frac{\d\oy}{\d\xi} = \oA(\xi;\la)\,\oy \quad\text{with}\quad
\oA = (B^{-1}AB - \mu_-^{[1]} I).
\end{equation}
The matrix $\oA(\xi;\la)$ in this equation is given by
\begin{subequations}
\label{oA}
\begin{equation}
\oA(\xi;\la) = \begin{bmatrix}
-\frac{1}{\ka} \phi_-(\xi) & -\frac{1}{\ka} \phi_-(\xi) \\[\jot]
\frac{1}{\ka} \phi_-(\xi) & -\ka+\frac{1}{\ka} \phi_-(\xi)
\end{bmatrix}, 
\end{equation}
where
\begin{equation}
\label{ka} 
\phi_-(\xi) = f'(\hu(\xi)) - f'(\hu_-) \text{ and }
\ka = \sqrt{c^2 + 4(\la-f'(\hu_-))}. 
\end{equation}
\end{subequations}
Note that the parameters~$c$ and~$\la$ are replaced by only one
parameter,~$\ka$. Now, suppose that~$\ou$ and~$\ov$ can be expanded in
inverse powers of~$\ka$:
\begin{align*}
\ou(\xi;\ka) &= \ou_0(\xi) + \ka^{-1} \ou_1(\xi) + \ka^{-2} \ou_2(\xi) +
\ka^{-3} \ou_3(\xi) + \cO(\ka^{-4}), \\
\ov(\xi;\ka) &= \ov_0(\xi) + \ka^{-1} \ov_1(\xi) + \ka^{-2} \ov_2(\xi) +
\ka^{-3} \ov_3(\xi) + \cO(\ka^{-4}). 
\end{align*}
If we substitute these expansions in~\eqref{bvpo} and equate the
coefficients of the powers of~$\ka$, we find:
\begin{align*}
0 &= 0, & 0 &= -\ov_0, \\ 
\ou_0' &= 0, & \ov_0' &= -\ov_1,  \\
\ou_1' &= -\phi_-(\xi) \, (\ou_0+\ov_0), &
\ov_1' &= \phi_-(\xi) \, (\ou_0+\ov_0) - \ov_2, \\
\ou_2' &= -\phi_-(\xi) \, (\ou_1+\ov_1), &
\ov_2' &= \phi_-(\xi) \, (\ou_1+\ov_1) - \ov_3.
\end{align*}
Assuming that $\ou(\xi)$ and $\ov(\xi)$ are bounded as
$\xi\to-\infty$, the solution of these equations (up to a
multiplicative constant) is
\begin{equation}
\label{oy-exp}
\begin{aligned}
\ou(\xi;\ka) &= 1 - \ka^{-1} \Phi_-(\xi) + \tfrac12 \ka^{-2}
\big(\Phi_-(\xi)\big)^2 + \cO(\ka^{-3}), \\
\ov(\xi;\ka) &= \ka^{-2} \phi_-(\xi) + \cO(\ka^{-3}), 
\end{aligned}
\end{equation}
where 
$$
\Phi_-(\xi) = \int_{-\infty}^\xi \phi_-(x) \,\d{x}. 
$$

We can do something similar to find the solution~$y$ of~\eqref{bvp3}
satisfying $y(\xi)\to0$ as $\xi\to+\infty$. Instead of~\eqref{oy}, we
write $y$ as
\begin{equation*}
y(\xi) = \exp(\mu_+^{[2]}\xi) \, B_+ \, \oy_+(\xi) \quad\text{where}\quad
B_+ = \begin{bmatrix} 1 & 1 \\ \mu_+^{[1]} & \mu_+^{[2]} \end{bmatrix}. 
\end{equation*}
Expanding $\oy_+$ in negative powers of~$\ka_+$, where 
\begin{equation}
\label{kaplus}
\ka_+ = \sqrt{c^2 + 4(\la-f'(\hu_+))}, 
\end{equation}
similar to~\eqref{ka}, we find that 
\begin{equation}
\label{oy-exp2}
\oy_+ = \begin{bmatrix} \ou_+ \\ \ov_+ \end{bmatrix} 
\text{ with }
\begin{cases}
\ou_+(\xi;\ka) = \ka_+^{-2} \phi_+(\xi) + \cO(\ka_+^{-3}), \\
\ov_+(\xi;\ka) = 1 + \ka_+^{-1} \Phi_+(\xi) + \tfrac12 \ka_+^{-2}
\big(\Phi_+(\xi)\big)^2 + \cO(\ka_+^{-3}), 
\end{cases}
\end{equation}
where
$$
\phi_+(\xi) = f'(\hu(\xi)) - f'(\hu_+)
\quad\text{and}\quad
\Phi_+(\xi) = \int_\xi^\infty \phi_+(x) \,\d{x}. 
$$
The Evans function is obtained by evaluating both the solution
satisfying $y(\xi)\to0$ as $\xi\to-\infty$ and the one satisfying
$y(\xi)\to0$ as $\xi\to+\infty$ at $\xi=0$, collecting the resulting
vectors in a matrix and computing the determinant of this matrix. This
yields
\begin{equation}
\label{dla-pre}
\begin{aligned}
D(\la) &= \big( B\oy(0) \big) \wedge \big( B_+\oy_+(0) \big) \\
&= 
\begin{bmatrix} 1 & 1 \\ \frac12(\ka-c) & -\frac12(\ka+c) \end{bmatrix}
\begin{bmatrix} \ou(0) \\ \ov(0) \end{bmatrix} \wedge
\begin{bmatrix} 1 & 1 \\ \frac12(\ka_+-c) & -\frac12(\ka_++c) \end{bmatrix}
\begin{bmatrix} \ou_+(0) \\ \ov_+(0) \end{bmatrix} \\
&= \tfrac12 (\ka-\ka_+) \big(\ov(0)\ov_+(0) - \ou(0)\ou_+(0) \big) \\
&\hspace{1cm} + \tfrac12 (\ka+\ka_+) \big(\ov(0)\ou_+(0) - \ou(0)\ov_+(0) \big)
\end{aligned}
\end{equation}
Substituting~\eqref{oy-exp} and~\eqref{oy-exp2} and using the fact
that $\ka-\ka_+ = \cO\big(|\la|^{-1/2}\big)$ as $|\la|\to\infty$, we find that
\begin{equation}
\label{dla-as}
\begin{aligned}
D(\la) &= -\tfrac12 (\ka+\ka_+) \ou(0)\ov_+(0) + \cO(\ka^{-2}) \\
&=  -2\la^{1/2} + \Phi - \tfrac14 \la^{-1/2} \Big(
\Phi^2 - 2f'(\hu_-) - 2f'(\hu_+) + c^2 \Big) + \cO(\la^{-1}), 
\end{aligned}
\end{equation}
where
$$
\Phi = \Phi_-(0) + \Phi_+(0) = \int_{-\infty}^0 \phi_-(x)
\,\d{x} + \int_0^\infty \phi_+(x) \,\d{x}. 
$$
The approach of Sandstede~\cite{sandstede:stability} yields $D(\la) =
-2\la^{1/2} + \cO(1)$. This agrees with~\eqref{dla-as}, but the
approach presented here gives two more terms. Furthermore, we can
easily find additional terms by extending the
expansions~\eqref{oy-exp}.

\section{Magnus methods}
\label{sec:magnus}

If we want to evaluate the Evans function numerically, we have to
solve the differential equation~\eqref{bvp}.
Moan~\cite{moan:efficient} studied methods based on the Magnus series
for the solution of Sturm--Liouville problems of the form $-(py')' +
qy = \la wy$ on a finite interval. Moan noticed that some of the
quantities involved in the computation are independent on the spectral
parameter~$\la$ and hence need to be computed only once when solving
the differential equation for several values of~$\la$. Moan also
proposed a modification of the method based on summing some of the
terms analytically which improves the accuracy when $\la$ is large.

J\'odar and Marletta~\cite{jodar.marletta:solving} noticed that the
Magnus method in combination with the compound matrix method performs
well on some scalar Sturm--Liouville problems of high order; see also
Greenberg and Marletta~\cite{greenberg.marletta:numerical}.

This approach was generalized by Aparicio, Malham and
Oliver~\cite{aparicio.malham.ea:numerical}, who proposed to use a
Magnus method for solving the boundary value problem~\eqref{bvp}. They
mentioned the robustness across all regimes as an advantage of Magnus
integrators. Furthermore, they pointed out that the computational cost
of Magnus methods, as well as other methods, can be decreased by means
of a precomputation technique, In this section, we further analyze the
behaviour of the Magnus method in the regime where $\la$ is large in
modulus.

Magnus~\cite{magnus:on} showed that the solution of the differential
equation $y' = A(\xi) y$ can be written as $y(\xi) = \exp(\Om(\xi)) \,
y(0)$, where the matrix $\Om(\xi)$ is given by the infinite series
\begin{equation}
\label{magnus}
\begin{aligned}
  \Om(\xi) &= \int_0^\xi A(x) \,\d x - \tfrac12 \int_0^\xi \left[
    \int_0^{x_1} A(x_2) \,\d x_2, A(x_1) \right] \d x_1 \\
  &\hspace{1cm} + \tfrac1{12} \int_0^\xi \left[
    \int_0^{x_1} A(x_2) \,\d x_2, \left[
      \int_0^{x_1} A(x_2) \,\d x_2, A(x_1)
    \right] \right] \d x_1 \\
  &\hspace{1cm} + \tfrac14 \int_0^\xi \left[
    \int_0^{x_1} \left[ \int_0^{x_2} A(x_3)
      \,\d x_3, A(x_2) \right] \d x_2, A(x_1) 
  \right] \d x_1 + \cdots,
\end{aligned}
\end{equation}
where $[\,\cdot\,,\,\cdot\,]$ denotes the matrix commutator defined by
$[X,Y] = XY - YX$. Moan and Niesen~\cite{moan.niesen:convergence}
proved that the series converges if $\int_0^\xi \|A(x)\| \,\d x <
\pi$.

The Magnus series can be used to solve linear differential equations
numerically, if we truncate the infinite series and approximate the
integrals numerically. For instance, if we retain only the first
term in the series and approximate $A(x)$ by the value at the
midpoint, we get $\Omega(\xi) \approx hA(\frac12\xi)$. The resulting
one-step method is defined by
\begin{equation}
\label{m2}
y_{k+1} = \exp\bigl(hA(\xi_k+\tfrac12h)\bigr) \, y_k,
\end{equation}
where $h$ denotes the step size and $y_k$ approximates the solution at
$\xi_k = \xi_0 + kh$. This method is called the \emph{Lie midpoint} or
\emph{exponential midpoint} method. It is a second-order method: the
difference between the numerical and the exact solution at a fixed
point~$\xi$ is~$\cO(h^2)$.

We can get a fourth-order method by truncating the Magnus
series~\eqref{magnus} after the second term. We replace the
matrix~$A(\xi)$ by the linear function $A_0 + \xi A_1$ which agrees
with~$A(\xi)$ at the two Gauss--Legendre points
\begin{equation}
\label{gl-points}
\xi_k^{[1]} = \xi_k + (\tfrac12-\tfrac16\sqrt3) h \quad\text{and}\quad
\xi_k^{[2]} = \xi_k + (\tfrac12+\tfrac16\sqrt3) h.
\end{equation}
This yields the scheme
\begin{subequations}
\label{m4}
\begin{equation}
y_{k+1} = \exp(\Om_k)\, y_k, 
\end{equation}
where
\begin{equation}
\Om_k = \tfrac12h \big( A(\xi_k^{[1]}) + A(\xi_k^{[2]}) \big) -
\tfrac{\sqrt3}{12} h^2 \big[ A(\xi_k^{[1]}), A(\xi_k^{[2]}) \big].
\end{equation}
\end{subequations}
The reader is referred to the review paper by Iserles, Munthe--Kaas,
N{\o}rsett and Zanna~\cite{iserles.munthe-kaas.ea:lie-group} for more
information on Magnus and related methods. 


If we define $\oy_k$ by $y_k = \exp(\mu_-^{[1]}\xi) \, B \, \oy_k$, as
suggested by~\eqref{oy}, then the fourth-order method~\eqref{m4}
transforms to
\begin{subequations}
\label{om4}
\begin{equation}
\label{om4a}
\oy_{k+1} = \exp(\oOm_k)\, \oy_k,
\end{equation}
where
\begin{equation}
\label{om4b}
\oOm_k = \tfrac12h \big( \oA(\xi_k^{[1]}) + \oA(\xi_k^{[2]}) \big) -
\tfrac{\sqrt3}{12} h^2 \big[ \oA(\xi_k^{[1]}), \oA(\xi_k^{[2]}) \big],
\end{equation}
\end{subequations}
with $\oA$ as given in~\eqref{oA}. So applying the Magnus method to
the transformed equation~\eqref{bvpo} and transforming the result back
to the original coordinate system is the same as applying it to the
original equation. This can be explained by the equivariance of the
Magnus method under linear transformations and exponential
rescalings~\cite{eng:on}.

Substitution of~\eqref{oA} in~\eqref{om4b} yields 
\begin{subequations}
\label{oOm+al+be}
\begin{equation}
\label{oOm}
\oOm_k = h \begin{bmatrix} 
  -\ka^{-1}\al_k & \be_k - \ka^{-1}\al_k \\
  \be_k + \ka^{-1}\al_k & -\ka+\ka^{-1}\al_k 
\end{bmatrix}
\end{equation}
with
\begin{equation}
  \label{al+be}
  \al_k = \tfrac12 \bigl( \phi_-(\xi_k^{[1]})+\phi_-(\xi_k^{[2]}) \bigr)
  \quad\text{and}\quad
  \be_k = -\frac{\sqrt{3}}{12} h \bigl(
  \phi_-(\xi_k^{[1]})-\phi_-(\xi_k^{[2]}) \bigr). 
\end{equation}
\end{subequations}
Note that $\al_k$ and $\be_k$ approximate $\phi_-$
and~$\frac1{12}h^2\phi'_-$ respectively. Furthermore, the exponential
midpoint rule~\eqref{m2} is also given by~\eqref{oOm}, but with $\al_k
= \phi_-(\xi_k+\tfrac12h)$ and $\be_k = 0$ instead of~\eqref{al+be}.

\subsection{Estimates for the local error}

The local error of a one-step method is the difference between the
numerical solution and the exact solution after one step. For the
Magnus method, the local error is 
$$
L_k = \exp(\Om_k) \, y(\xi_k) - y(\xi_{k+1}), 
$$
or, in transformed coordinates,
\begin{equation}
\label{oLk-dfn}
\oL_k = \exp(\oOm_k) \, \oy(\xi_k) - \oy(\xi_{k+1}).
\end{equation}
The exponential of the matrix~$\oOm_k$ is most easily calculated by
diagonalization: if $\oOm_k = V_k \La_k V_k^{-1}$ with
$\La_k$~diagonal, then $\exp(\oOm_k) = V_k \, \exp(\La_k) \, V_k^{-1}$
and $\exp(\La_k)$ is formed by simply exponentiating the entries on
the diagonal. In fact, the diagonal entries of~$\La_k$ are
\begin{equation}
\label{lak}
\begin{aligned}
\la_k^{[1]} &= h \Bigl( \ka^{-1} \big(\be_k^2-\al_k\big) -
\ka^{-3} \big(\be_k^2-\al_k\big)^2 + \cO(\ka^{-5}) \Bigr), \\
\la_k^{[2]} &= -h \Bigl( \ka + \ka^{-1} \big(\be_k^2-\al_k\big) -
\ka^{-3} \big(\be_k^2-\al_k\big)^2 + \cO(\ka^{-5}) \Bigr), 
\end{aligned}
\end{equation}
as $|\ka| \to \infty$. The definition of~$\ka$ in~\eqref{ka} implies
that $\Re\ka > 0$ unless $\la$ is real and $\la \le f'(\hu_-) -
\tfrac14 c^2$. Under this condition, $-\Re \la_k^{[2]} \gg 1$ if $h|\ka|
\gg 1$ and hence $\exp(\la_k^{[2]})$ is exponentially small.

We now assume that we are in the regime with $|\la| \gg h^{-2}$, $h
\to 0$, and $\la$ bounded away from the negative real axis in the
sense that $|\arg \la| < \pi-\eps$ where $\eps > 0$. In this regime,
$\exp(\la_k^{[2]})$ is exponentially small. Taking this into account,
a lengthy but straightforward calculation shows that
\begin{equation}
\label{eoOm}
\exp(\oOm_k) = \begin{bmatrix}
  1 - \dfrac{h\chi_k}{\ka} + \cO(\ka^{-2})
  & \dfrac{\be_k}\ka - \dfrac{\al_k+h\be_k\chi_k}{\ka^2} +
  \cO(\ka^{-3}) \\[2\jot]
  \dfrac{\be_k}\ka + \dfrac{\al_k+h\be_k\chi_k}{\ka^2} +
  \cO(\ka^{-3}) 
  & \dfrac{\be_k^2}{\ka^2} - \dfrac{h\be_k^2\chi_k}{\ka^3} +
  \cO(\ka^{-4}) 
\end{bmatrix},
\end{equation}
where $\chi_k = \al_k-\be_k^2$.  Substituting this result and the
approximation~\eqref{oy-exp} for the exact solution in the
definition~\eqref{oLk-dfn}, and using the definitions of~$\chi_k$,
$\beta_k$, and $\Phi_-$, we find that the local error of the Magnus
method is given by
\begin{equation}
\label{oLk}
\oL_k = \begin{bmatrix} 
\ka^{-1} \ga_k + \cO(\ka^{-2}h^4) \\[\jot]
\ka^{-1} \be_k + \cO(\ka^{-2}h) 
\end{bmatrix}, 
\end{equation}
where
\begin{equation}
\label{ga}
\begin{aligned}
\ga_k &= 
\int_{\xi_k}^{\xi_{k+1}} \phi_-(x) \,\d{x} - h(\al_k-\be_k^2) \\
&= h^5 \left( \tfrac{1}{4320} \phi''''(\xi_k+\tfrac12h) + 
\tfrac1{144} \big(\phi'(\xi_k+\tfrac12h)\big)^2 \right) + \cO(h^7) 
\end{aligned}
\end{equation}
and
\begin{equation}
\label{be}
\be_k = \tfrac1{12} h^2 \phi'(\xi_k+\tfrac12h) + \cO(h^4).
\end{equation}
In deriving the above expression, we also replaced $\phi_-(\xi)$ by $\phi(\xi)
= f'(\hu(\xi))$. This is allowed since they differ by a constant term,
and only the derivative appears in~\eqref{oLk}. 

Equation~\eqref{oLk} gives the local error of the fourth-order Magnus
method~\eqref{m4} in transformed variables. Since the method has order
four, the local error is~$\cO(h^5)$ as $h\to0$ when solving a fixed
equation.  However, in our case, the constraint $|\la| \gg h^{-2}$
implies that $\la$, and hence the relative influence of the
coefficients in the equation, must change as $h$ approaches zero. It
turns out that the local error is~$\cO(h^2)$ in this setting.  In
other words, the method behaves like a first-order method (globally).
This phenomenon is called \emph{order reduction}.

The cause of this order reduction is the stiffness of the differential
equation~\eqref{bvp3}. Indeed, if we define the \emph{stiffness ratio}
as the quotient between the largest and smallest
eigenvalue (as in Iserles~\cite{iserles:first}), then the stiffness
quotient is 
$$ 
\left| \frac{\la_k^{[2]}}{\la_k^{[1]}} \right| =
\left| \frac{\ka^2}{\be_k^2-\al_k} \right| + \cO(1), 
$$
where the leading term shown grows like~$|\la|$. Hence, the problem is
stiff if $\la$ is large in modulus, causing troubles for the numerical
method.

\subsection{Estimates for the global error}
\label{sec:globerr}

The global error is the error of the numerical method after
\emph{several} steps, say~$k$. Hence, the global error is $E_k = y_k -
y(\xi_k)$, with $y_k$ defined by the numerical method starting from
$y_0 = y(\xi_0)$. For the Magnus method~\eqref{m4}, the global error
satisfies the recursion relation $E_{k+1} = \exp(\Om_k) \, E_k + L_k$
with $E_0 = 0$, or, in transformed coordinates,
\begin{equation}
\label{ge-rec}
\oE_{k+1} = \exp(\oOm_k) \, \oE_k + \oL_k, \quad \oE_0 = 0.
\end{equation}
A routine induction argument using~\eqref{eoOm} and~\eqref{oLk} shows
that the leading term of the global error is given by
\begin{equation}
\label{oEk-ih}
\oE_k = \begin{bmatrix} 
  \ka^{-1} \sum_{j=0}^{k-1} \ga_j + \cO(\ka^{-2}h^4) \\[\jot] 
  \ka^{-1} \be_{k-1} + \cO(\ka^{-2}h) 
\end{bmatrix}.
\end{equation}
This shows an advantage of stiffness: the exact flow quickly reduces
the error in the stiff component. The Magnus method inherits this
property here and annihilates at every step the error in the stiff
component up to leading order (if $|\ka| h \gg 1$). On the other hand,
the first (nonstiff) component of the error is propagated without
change. Since the local error in the stiff component is much bigger
than the error in the nonstiff component (order~$h^2$ versus
order~$h^5$), we arrive at the surprising conclusion that the local
and global error are equal at leading order.

Substituting the definitions of~$\al_k$, $\be_k$, and $\ga_k$
in~\eqref{al+be} and~\eqref{ga} and approximating the sum by an
integral, we find that
\begin{equation}
\label{oEk}
\oE_k = \begin{bmatrix}
  \ka^{-1}h^4 \Bigl( \frac1{4320} \bigl(\phi'''(\xi_k)-\phi'''(\xi_0)\bigr) 
  + \frac1{144} \int_{\xi_0}^{\xi_k} (\phi'(\xi))^2 \,\d{\xi} \Bigr) +
  \cO(\ka^{-1}h^6,\ka^{-2}h^4) \\[5pt] 
  \frac1{12} \ka^{-1}h^2 \phi'(\xi_k-\frac12h) +
  \cO(\ka^{-1}h^4,\ka^{-2}h)
\end{bmatrix}.  
\end{equation}
We see that the global error is~$\cO(h^2)$. Usually, a factor~$h$ is
lost in the transition from the local to the global error, but here
both the local and the global error are~$\cO(h^2)$, because the local
error is mainly in the stiff component. Hence, the fourth-order Magnus
method given in~\eqref{m4} behaves like a second-order method if one
considers the global error.  The numerical experiments in
Section~\ref{sec:num} support this analysis. However, the fact that
the global error is~$\cO(h^2)$ does not tell the whole story, as the
relative error (the error divided by the magnitude of the solution)
may give a better picture. Indeed, it follows from~\eqref{oy}
and~\eqref{oy-exp} that the solution~$y$ grows as~$\sqrt{\la}$, so the
relative error is $\cO(h^2/\sqrt{\la})$. In other words, the relative
error decreases as $|\la|\to\infty$.

We can compute the error associated with the solution satisfying the
boundary condition that $y(\xi)\to0$ as $\xi\to+\infty$
similarly. Instead of~\eqref{oEk-ih}, we now have
\begin{equation}
\label{oEkp-ih}
\oE^+_k = \begin{bmatrix} 
  \ka_+^{-1} \sum_{j=0}^{k-1} \ga^+_j + \cO(\ka_+^{-2}h^4) \\[1.5\jot] 
  \ka_+^{-1} \be^+_{k-1} + \cO(\ka_+^{-2}h) 
\end{bmatrix},
\end{equation}
where $\al_k^+$, $\be_k^+$ and $\ga_k^+$ are as given
in~\eqref{al+be}, with $\phi_+$ and~$\Phi_+$ replacing~$\phi_-$
and~$\Phi_-$, respectively, and $\ka_+$ is as given
in~\eqref{kaplus}.

\subsection{The error in the Evans function}
\label{sec:everr}

The Evans function given in~\eqref{dla-pre} is
$$
D(\la) = \big( B\oy(0) \big) \wedge \big( B_+\oy_+(0) \big).
$$
The numerical error when evaluating the Evans function is therefore
\begin{equation}
\label{ED1}
E_D = \big( B\oy(0) \big) \wedge \big( B_+\oE_k^+ \big) 
+ \big( B\oE_k \big) \wedge \big( B_+\oy_+(0) \big) 
+ \big( B\oE_k \big) \wedge \big( B_+\oE_k^+ \big).
\end{equation}
We can expand this in the same manner as in~\eqref{dla-pre}. The first
term on the right-hand side becomes
$$
\tfrac12 (\ka-\ka_+) \big([\oE]_2\ov_+(0) - [\oE]_1\ou_+(0) \big)
+ \tfrac12 (\ka+\ka_+) \big([\oE]_2\ou_+(0) - [\oE]_1\ov_+(0) \big).
$$
The dominating term in this expression is $-\tfrac12 (\ka+\ka_+)
[\oE]_1\ov_+(0)$, which is~$\cO(\la^0h^4)$; all other terms
are~$\cO(\la^{-1}h^2)$ or smaller (recall that we assumed that $|\la|
\gg h^{-2}$). Therefore, the dominating contribution to the error in
the Evans function comes from the first (nonstiff) component of the
global error, even though the second (stiff) component is larger. This
is because the stiff and nonstiff directions are exchanged when you
integrate in the other direction.  Hence, when taking the wedge
product of the global error of the solution on $[-\infty,0]$ with the
solution itself on $[0,+\infty]$, the stiff component of the global
error is paired with the stiff component of the solution; similarly,
the nonstiff component of the global error is paired with the nonstiff
component of the solution.  Since the solution is mainly along the
nonstiff direction, the nonstiff component of the global error (which
has order~$h^4$) is brought to the fore, and the stiff component of
the global error (which has order~$h^2$) is reduced.

Substituting the exact solution from~\eqref{oy-exp}
and~\eqref{oy-exp2} and the global error from~\eqref{oEk-ih}
and~\eqref{oEkp-ih} in~\eqref{ED1}, we find that the error in the
Evans function is
$$
E_D = \sum_{j=0}^{k-1} \ga_j + \sum_{j=0}^{k-1} \ga^+_j +
\cO(\la^{-1/2}h^4).
$$ 
Assuming that the differential equation is solved on the
intervals~$[-L,0]$ and~$[0,L]$, with $L = Nh$, this evaluates to
\begin{multline}
\label{ED2}
E_D = \left[ h \sum_{j=-N}^{N-1} \Big(
  \phi\big(jh+(\tfrac12-\tfrac16\sqrt3\big)h +
  \phi\big(jh+(\tfrac12+\tfrac16\sqrt3\big)h \Big) 
- \int_{-L}^L \phi(x) \,\d{x} \right] \\
- h \sum_{j=0}^{k-1} \be^2_j - h \sum_{j=0}^{k-1} (\be^+_j)^2 +
\cO(\la^{-1/2}h^4).
\end{multline}
The term within brackets is difference between the approximation of
$\int_{-L}^L \phi(x) \,\d{x}$ by two-point Gauss--Legendre quadrature
and the integral itself. In our setting, all derivatives of the
travelling wave~$\hu$, and therefore also of the function~$\phi = f'
\circ \hu$, vanish at infinity (see~\S\ref{sec:evans}). Now, assuming
that $L$ is so large that the derivatives of~$\phi$ at~$L$ are
negligible, the error in Gauss--Legendre quadrature vanishes at all
orders in~$h$, for essentially the same reason that the trapezoidal
rule is so effective for periodic integrands; this is easily proved
with the Euler--MacLaurin formula (see, for instance, Davis and
Rabinowitz~\cite[\S3.4]{davis.rabinowitz:numerical}). Hence, only the
sums involving the $\be_j$ and $\be^+_j$ survive. These can be
approximated easily using~\eqref{be}, and we find that
\begin{equation}
\label{ED}
E_D = -\frac{h^4}{144} \int_{-\infty}^\infty \big(\phi'(x)\big)^2
\,\d{x} + \cO(h^6)
\quad\text{with}\quad \phi(\xi) = f'(\hu(\xi)).
\end{equation}
So, in the end, the error in the Evans function is of order~$h^4$,
which is just what one would expect from a fourth-order method.

\subsection{The exponential midpoint rule}
\label{sec:emp}

We saw above that the fourth-order Magnus method~\eqref{m4} suffers
severe order reduction when solving~\eqref{bvp3} with $|\la| \gg
1$. This is not the case for all methods. There are even methods based
on the Magnus series which do not suffer global order reduction, like
the exponential midpoint rule~\eqref{m2} which has order two (this
method is also known as the second-order Magnus method). When applied
to~\eqref{bvp3}, this method is of the form \eqref{om4a},~\eqref{oOm}
with $\al_k$ and~$\be_k$ given by $\al_k = \phi_-(\xi_k+\tfrac12h)$
and $\be_k = 0$ respectively. If we substitute this in~\eqref{oLk}, we
see that the $\ka^{-1}$ term in the second (stiff) component drops
out, and that the local error is
$[\,\cO(\ka^{-1}h^3)\,\,\cO(\ka^{-2}h)\,]^\top$. So, the exponential
midpoint rule does suffer some local order reduction, but not as
severe as the fourth-order Magnus method, for which the second
component of the local error is of order~$\ka^{-1}h^2$.

The global error can be computed as in~\S\ref{sec:globerr}.  Again,
only the nonstiff component propagates, so the global error is
$[\,\cO(\ka^{-1}h^2)\,\,\cO(\ka^{-2}h)\,]^\top$. As $|\ka| \gg
h^{-1}$, the first component dominates and the exponential midpoint
rule effectively does not suffer from order reduction if one looks at
the global error.

Continuing to find the error in the Evans function, as we did
in~\S\ref{sec:everr}, we find
\begin{equation}
\label{ED-m2}
E_D = h \sum_{j=-N}^{N-1} \phi\big(jh+\tfrac12h)
- \int_{-L}^L \phi(x) \,\d{x} + \cO(\la^{-1/2}h^2).
\end{equation}
Comparing with~\eqref{ED2} for the fourth-order Magnus method, we see
that the sums involving the~$\be_j$ and~$\be^+_j$ have dropped out
(because $\be_j=0$), and that the two-point Gauss--Legendre quadrature
is replaced by the trapezoidal rule. Again, the error in the
trapezoidal rule vanishes at all orders if $L$ is sufficiently
large. Hence, the error in the Evans function
is~$\cO(\la^{-1/2}h^2)$. In contrast, the fourth-order Magnus method
has $E_D = \cO(h^4)$, see~\eqref{ED}. Thus, we can expect the
second-order method to be more accurate than the fourth-order
method. The experiments in Section~\ref{sec:num} confirm this.

\section{The Gauss--Legendre method}
\label{sec:gl}

Most numerical computations of the Evans function reported in the
literature use a Runge--Kutta method, in particular the classical
explicit fourth-order method and the two-stage Gauss--Legendre
method. As the differential equation that we want to solve, is stiff,
we consider the two-stage Gauss--Legendre method. The method is given
by
\begin{equation}
\label{gl4}
\begin{aligned}
s_1 &= A(\xi_k^{[1]}) \, 
\bigl( y_k + \tfrac14h s_1 + (\tfrac14-\tfrac{\sqrt3}6)h s_2 \bigr), \\
s_2 &= A(\xi_k^{[2]}) \, 
\bigl( y_k + (\tfrac14+\tfrac{\sqrt3}6)h s_1 + \tfrac14h s_2 \bigr), \\
y_{k+1} &= y_k + \tfrac12h (s_1 + s_2),
\end{aligned}
\end{equation}
where $\xi_k^{[1]}$ and $\xi_k^{[2]}$ are the Gauss--Legendre points,
given in~\eqref{gl-points}. 

We will now analyse the error committed by the Gauss--Legendre method
when computing the Evans function. After the usual coordinate
transformation, cf.~\eqref{oy}, and substitution of the matrix given
in~\eqref{oA}, we can solve the system~\eqref{gl4}. After a lengthy
but relatively straightforward calculation, we find that $\oy_{k+1} =
\Psi_k \oy_k$ with 
$$
\Psi_k = \begin{bmatrix}
1 - \ka^{-1}h\al_k + \frac12\ka^{-2}h^2\al_k^2 
& 12 \ka^{-2}h^{-1}\be_k \\[2\jot]
12 \ka^{-2}h^{-1}\be_k 
& 1 - 12\ka^{-1}h^{-1} + 72\ka^{-2}h^{-2}
\end{bmatrix} + \cO(\ka^{-3}) 
$$
where $\al_k$ and~$\be_k$ are as defined in~\eqref{al+be}.  The local
error can now be found by substituting this matrix and the exact
solution in $\oL_k = \Psi_k\oy(\xi_k) - \oy(\xi_{k+1})$,
cf.~\eqref{oLk-dfn}. This yields
\begin{equation}
\label{oLk-gl4}
\oL_k = \begin{bmatrix} 
  \ka^{-1} \oL_k^a + \ka^{-2} \oL_k^b + \cO(\ka^{-3}h^5) \\[\jot]
  \ka^{-2} \oL_k^c + \cO(\ka^{-3}h^2) 
\end{bmatrix}
\end{equation}
where
\begin{subequations}
\begin{align}
\label{oLka-gl4}
\oL_k^a &= \int_{\xi_k}^{\xi_k+h} \phi_-(x) \,\d{x}  - h\al_k = \cO(h^5) \\
\label{oLkb-gl4}
\oL_k^b &= \Phi_-(\xi_k) \bigg( h\al_k - \int_{\xi_k}^{\xi_k+h} \!
\phi_-(x) \,\d{x} \biggr) 
+ \tfrac12(h\al_k)^2 - \tfrac12\bigg(\int_{\xi_k}^{\xi_k+h} \!
\phi_-(x) \,\d{x} \bigg)^2  \\
\notag
&= - \oL_k^a \big( \Phi_-(\xi_k) + h\al_k + \tfrac12\oL_k^a \big) =
\cO(h^5) \\[2mm]  
\label{oLkc-gl4}
\oL_k^c &= 12h^{-1}\be_k - \phi_-(\xi_k+h) + \phi_-(\xi_k) = \cO(h^3).
\end{align}
\end{subequations}
For reasons which will soon become clear, we must retain the
$\ka^{-2}$ term in the above expression, in contrast to~\eqref{oLk}
for the Magnus method. We see that the two terms in the first
(nonstiff) component are of order~$\ka^{-1}h^5$ and~$\ka^{-2}h^5$,
while the second (stiff) component is of order~$\ka^{-2}h^3$. This is
a similar situation as with the exponential midpoint rule, except that
\eqref{gl4} is a fourth-order method.

To find the global error, we solve the recursion relation $\oE_{k+1} =
\Psi_k\oE_k + \oL_k$, $\oE_0=0$, cf.~\eqref{ge-rec}. The solution is
\begin{equation}
\label{oEk-gl4}
\oE_k = \begin{bmatrix}
\ka^{-1} \sum_{j=0}^{k-1} \oL_j^a + \ka^{-2} \sum_{j=0}^{k-1} \big(
\oL_j^b - h \al_j \sum_{i=0}^{j-1} \oL_i^a \big) + \cO(\ka^{-3}h^4) \\[2\jot]
\ka^{-2} \sum_{j=0}^{k-1} \oL_j^c + \cO(\ka^{-3}h^2) 
\end{bmatrix}
\end{equation}
Again, only the error in the nonstiff component propagates.

Finally, we compute the error in the Evans function. Estimating the
various terms in~\eqref{ED1}, we find that
\begin{equation}
\label{ED1-gl4}
E_D = -\tfrac12 (\ka+\ka_+) \big([\oE]_1\ov_+(0) + \ou(0)\,[\oE^+]_2
\big) + \cO(\la^{-1/2}h^8,\la^{-3/2}h^2).
\end{equation}
Hence, we wish to compute $\cX = [\oE]_1\ov_+(0) + \ou(0)\,[\oE^+]_2$.
Substitution of the exact solution, given in~\eqref{oy-exp}
and~\eqref{oy-exp2}, and the global error~\eqref{oEk-gl4} yields
$$
\cX
= \ka^{-1} \sum_{j=0}^N \big( \oL_j^a + \oL_j^{a,+} \big) 
+ \ka^{-2} \cX_2 + \cO(\ka^{-3}h^4). 
$$
where
$$
\cX_2 = \sum_{j=0}^N \bigg( \oL_j^b - h\al_j \sum_{i=0}^{j-1}
\oL_i^a - \oL_j^a \, \Phi_+(0) 
+ \oL_j^{b,+} - h\al_j^+ \sum_{i=0}^{j-1} \oL_i^{a,+} -
\oL_j^{a,+}\,\Phi_-(0) \bigg).
$$
The sum $\sum_{j=0}^N \big( \oL_j^a + \oL_j^{a,+} \big)$ is the term
within brackets in~\eqref{ED2}, so again, it vanishes at all orders
in~$h$. For the $\Phi_-$ and $\Phi_+$ terms, we use
\begin{equation}
\label{Phim-approx}
\Phi_-(\xi_j) = \int_{-\infty}^{\xi_j} \phi_-(x) \,\d{x} \approx
\sum_{i=0}^{j-1} \int_{\xi_i}^{\xi_i+h} \phi_-(x) \,\d{x} = 
\sum_{i=0}^{j-1} \big( \oL_i^a + h\al_i \big),
\end{equation}
where the approximate equality becomes exact in the limit
$L\to\infty$; again, we assume that $L$ is so large that we can
neglect any errors here. Substitution of~\eqref{Phim-approx}, its
analogue for~$\Phi_+$, and~\eqref{oLka-gl4} and~\eqref{oLkb-gl4} yields
\begin{multline*}
\cX_2 = - \sum_{j=0}^N \bigg( h\al_j\oL_j^a + h\al_j^+\oL_j^{a,+} +
\sum_{i=0}^N \big(\al_i^+ \oL_j^a + \al_i \oL_j^{a,+} \big) \\
+ \sum_{i=0}^{j-1} \big( h\al_i\oL_j^a + h\oL_i^a\al_j +
h\al_i^+\oL_j^{a,+} + h\oL_i^{a,+}\al_j^+ \big) \bigg) + \cO(h^8),
\end{multline*}
where the remainder term comes from estimating terms like $\sum_j
\sum_i \oL_i^a \oL_j^a$. This nested sum can be written as the product
of two sums:
$$
\cX_2 = \sum_{i=0}^N h(\al_i + \al_i^+) \cdot \sum_{j=0}^N \big(
\oL_j^a + \oL_j^{a,+} \big) + \cO(h^8).
$$
So we arrive again at the sum $\sum_{j=0}^N \big( \oL_j^a +
\oL_j^{a,+} \big)$, which vanishes at all orders. 

Substituting everything back into~\eqref{ED1-gl4}, we find that the
error in the Evans function is given by
\begin{equation}
\label{ED-gl4}
E_D = \cO(\la^{-1/2}h^8, \la^{-1}h^4, \la^{-3/2}h^2).
\end{equation}
This is clearly better than the fourth-order Magnus method, with $E_D
= \cO(\la^0h^4)$, and the exponential midpoint rule, with $E_D =
\cO(\la^{-1/2}h^2)$.

\section{Numerical experiments}
\label{sec:num}

In this section, we evaluate the Evans function for a particular
example. The error in this computation is determined and compared
against the estimates derived in the previous sections. 

The example is the Fisher equation
\begin{equation}
\label{fisher}
u_t=u_{xx}+u-u^2.
\end{equation}
This is a reaction--diffusion equation of the form~\eqref{rd}.
Fisher~\cite{fisher:wave} used it to describe the transmission of
genes in a population.  It is now viewed as the prototype equation
admitting travelling front
solutions~\cite[\S11.2]{murray:mathematical}.

The Fisher equation supports a travelling wave solution with wave
speed $c = -\frac56\sqrt6$. In fact, this solution is known
analytically:
$$
u(x,t) = \hu(\xi) = \frac{1}{\big(1+\e^{\xi/\sqrt6}\big)^2}
\quad\text{where}\quad \xi = x - \tfrac56\sqrt6 \, t.
$$
Suppose that we wish to determine the stability of this travelling
wave. We are led to consider the eigenvalue problem~\eqref{bvp3},
which in this case reads
\begin{subequations}
\label{fisher-bvp}
\begin{equation}
\frac{\d y}{\d\xi} = \begin{bmatrix}
0 & 1 \\ \la - \phi(\xi) & -c
\end{bmatrix} y, 
\end{equation}
where
\begin{equation}
\phi(\xi) = 1-2\hu(\xi) = 1 - \frac{2}{\big(1+\e^{\xi/\sqrt6}\big)^2}.
\end{equation}
\end{subequations}
We solve this equation with the fourth-order Magnus method, given
by~\eqref{m4}. In the previous section, we derived the local error
estimate~\eqref{oLk}. For the Fisher equation, this estimate evaluates
to
\begin{equation}
\label{fisher-oLk}
\oL_k \approx \frac1\ka \begin{bmatrix}
\dfrac{h^5 \e^{\xi/\sqrt6} \bigl( -8\e^{3\xi/\sqrt6} + 33 \e^{2\xi/\sqrt6} +
  702 \e^{\xi/\sqrt6} + 1 \bigr)}%
{38880 \bigl( 1 + \e^{\xi/\sqrt6} \bigr)^6} \\[5mm]
\dfrac{\sqrt6 \, h^2 \e^{\xi/\sqrt6}}{18 \bigl( 1 + \e^{\xi/\sqrt6} \bigr)^3}
\end{bmatrix}, 
\end{equation}
where $\ka = \sqrt{c^2+4(\la+1)}$ and $\xi$ is short for~$\xi_k$.

The global error estimate is given in~\eqref{oEk}. Assuming that
$\xi_0$ is negative and so large in magnitude that we can take $\xi_0
= -\infty$, we find that
\begin{equation}
\label{fisher-oEk}
\oE_k \approx \frac1\ka \begin{bmatrix}
\dfrac{\sqrt{6} \, h^4 \e^{\xi/\sqrt6} \bigl( 36\e^{4\xi/\sqrt6} +
  180\e^{3\xi/\sqrt6} + 364 \e^{2\xi/\sqrt6} + 353 \e^{\xi/\sqrt6} + 1 \bigr)}%
{38880 \bigl( 1 + \e^{\xi/\sqrt6} \bigr)^6} \\[5mm]
\dfrac{\sqrt6 \, h^2 \e^{\xi/\sqrt6}}{18 \bigl( 1 + \e^{\xi/\sqrt6} \bigr)^3}
\end{bmatrix}.
\end{equation}
Finally, estimate~\eqref{ED} for the error in the Evans function is
\begin{equation}
\label{fisher-ED}
E_D \approx -\frac{\sqrt6}{1080} h^4 \approx -0.002268 \, h^4.
\end{equation}
This estimate is independent of the parameter~$\ka$.

\begin{figure}
  \begin{center}
    \includegraphics[width=1.0\linewidth]{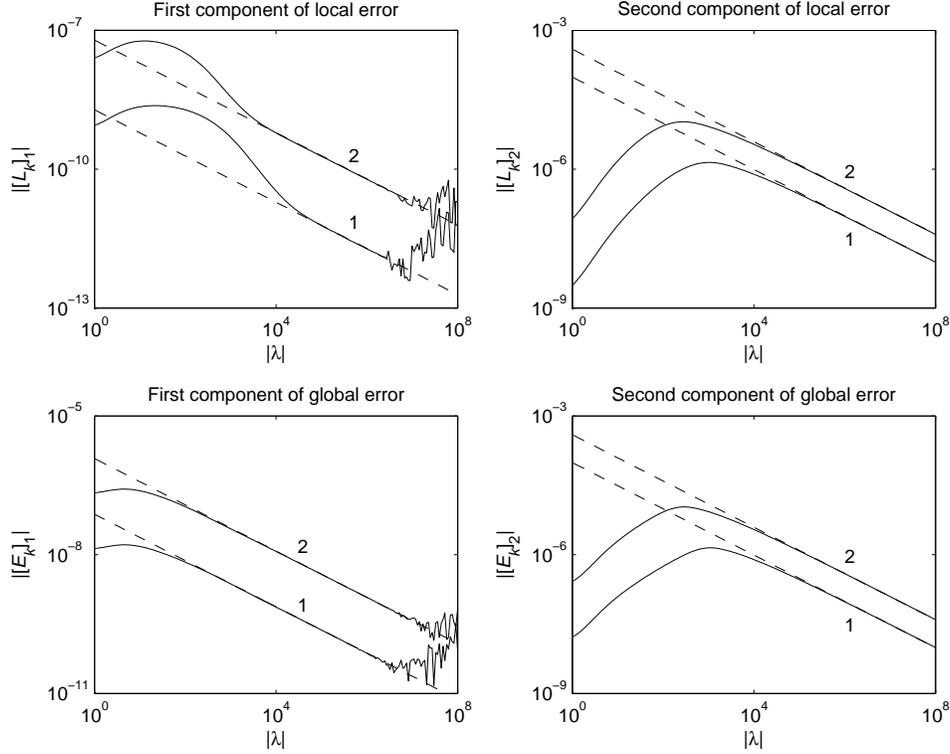}
  \end{center}
  \caption{The solid lines in the graphs on the top row show the local
    error committed by the fourth-order Magnus method. The step
    size~$h$ is 0.1 and~0.2 for the line labelled~1 and~2,
    respectively. The dash lines show the local error
    estimate~\eqref{fisher-oLk}. On the bottom row, the solid lines
    shows the global error and the dash lines show the
    estimate~\eqref{fisher-oEk}.}
  \label{scalar_fig1}
\end{figure}

\begin{figure}
  \begin{center}
    \includegraphics[width=1.0\linewidth]{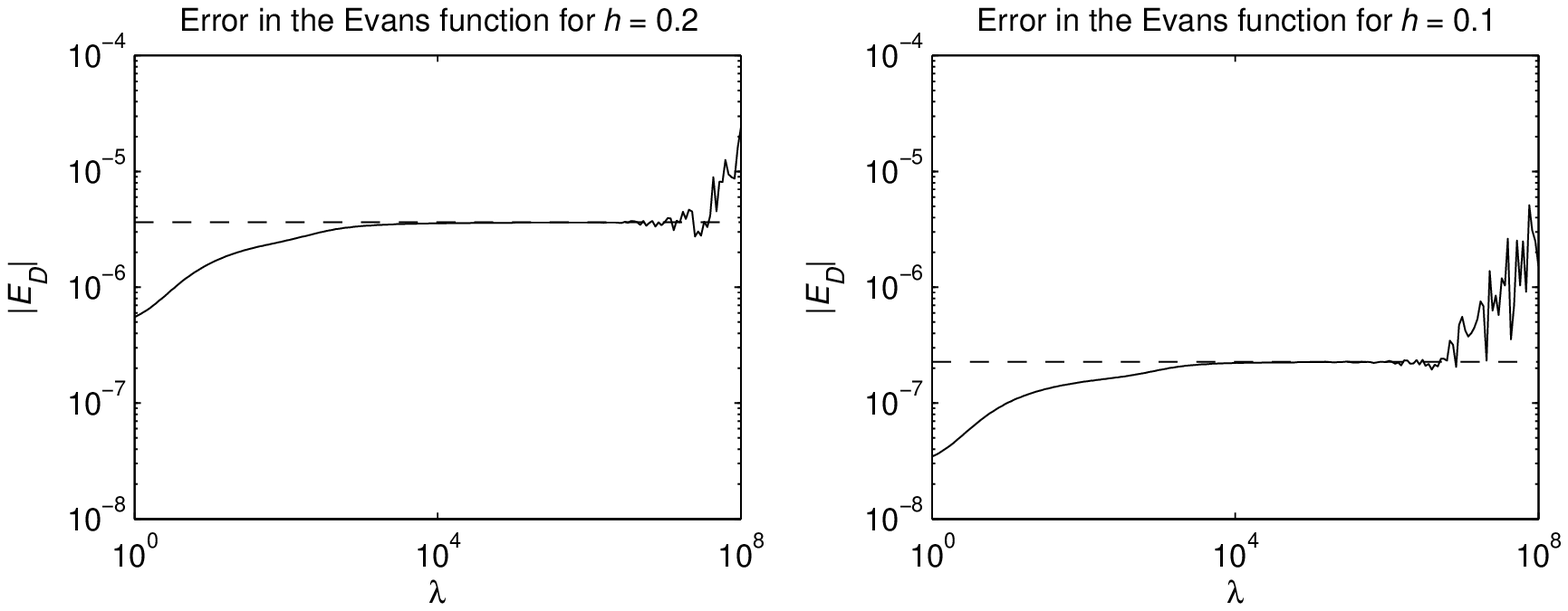}
  \end{center}
  \caption{The solid line shows the error in the Evans function as
    evaluated by the fourth-order Magnus method, while the dash line
    shows the error estimate~\eqref{fisher-ED}. The step size is
    $h=0.2$ for the graph on the left and $h=0.1$ for the graph on the
    right.}
  \label{scalar_fig2}
\end{figure}

We perform some numerical experiments to check the validity of these
estimates. First, we solve~\eqref{fisher-bvp} from $\xi=-30$ using the
fourth-order Gauss--Legendre method with step size $h = 0.02$. We will
refer to this solution as the ``exact'' solution. Then, we take the
``exact'' solution at $\xi=-1$ and do a single step with the
fourth-order Magnus method with step size $h=0.2$ or $h=0.1$. The
local error can now be determined by comparing the result of this
single step against the ``exact'' solution; this local error is
plotted in the top row of Figure~\ref{scalar_fig1}, together with the
local error estimate~\eqref{fisher-oLk}. The horizontal axis in the
plots shows the imaginary part of the eigenvalue parameter~$\la$,
which varies from~$\i$ to $10^8\,\i$ in our experiments.

The global error can be determined by solving~\eqref{fisher-bvp} from
$\xi=-30$ till $\xi=-1$ with the fourth-order Magnus method and
comparing it against the ``exact'' solution. This results in the
bottom row of Figure~\ref{scalar_fig1}. Finally,
Figure~\ref{scalar_fig2} shows the difference between the Evans
function as computed by the Magnus method and the ``exact'' value,
compared against the estimate~\eqref{fisher-ED}.

All graphs show that the error estimates agree well with the actual
error when $\la$ is moderately large in magnitude. However, the
numerical method starts to break down when $|\la|$ increases
above~$10^7$.

\label{expm}
The implementation used in the experiments is a straightforward
\textsc{Matlab} code. One detail proved to be important, namely, the
computation of the matrix exponential in~\eqref{m4}. The standard
routine for this is called \verb|expm| and uses Pad\'e approximation
combined with scaling and squaring. However, we found that an
alternative approach based on the Schur decomposition and implemented
in the \textsc{Matlab} routine \verb|expmdemo3| works better in our
case. Specifically, when using Pad\'e approximation, the numerical
method loses accuracy around $|\la| = 10^5$, as opposed to $|\la| =
10^7$ for the Schur decomposition. Generally, the Schur decomposition
runs into trouble when the matrix to be exponentiated is nearly
defective, but in our case the eigenvalues are far apart,
cf.~\eqref{lak}. The reader is refered to the article by Moler and Van
Loan~\cite{moler.loan:nineteen*1} for an extensive discussion on this
subject.

\begin{figure}
  \begin{center}
    \includegraphics[width=1.0\linewidth]{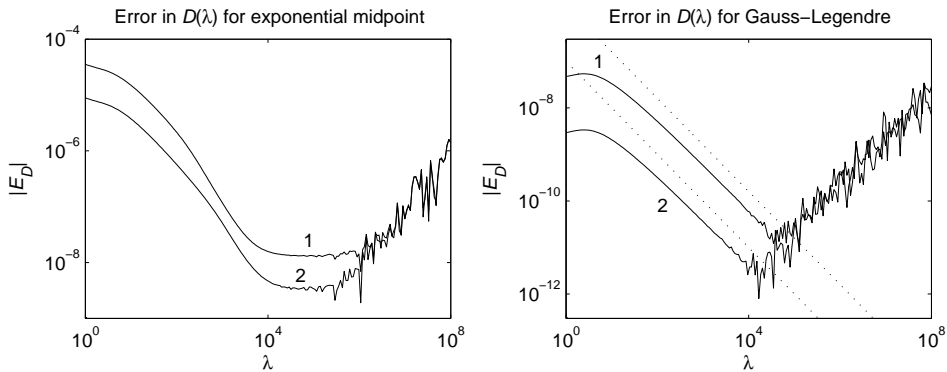}
  \end{center}
  \caption{The left graph shows the error in the Evans function as
  evaluated by the exponential midpoint rule~\eqref{m2}, while the
  right graph shows the same for the fourth-order Gauss--Legendre
  method~\eqref{gl4}. The step size is $h=0.2$ for the curve
  labelled~2 and $h=0.1$ for the curve labelled~1. The dotted lines in
  the second plot show $10^{-3}h^4/|\la|$.}
  \label{scalar_fig3}
\end{figure}

When the experiment is repeated with the exponential midpoint
rule~\eqref{m2} and the fourth-order Gauss--Legendre
method~\eqref{gl4}, the error in the Evans function is as plotted in
Figure~\ref{scalar_fig3}. We concluded in Section~\ref{sec:emp} that,
because the exponential midpoint rule suffers less from order
reduction, it is likely to have a smaller error than the fourth-order
Magnus method if $|\la|$ is large. The error plots confirm
this. 

For the fourth-order Gauss--Legendre method, which is the method that
is more relevant in practice, the error in the Evans function is shown
in the right half of Figure~\ref{scalar_fig3}. We found the
estimate~\eqref{ED-gl4} for the error, and the numerical results show
that the term of order $\la^{-1}h^4$ dominates: the line $10^{-3}
|\la|^{-1} h^4$ tracks the graphs closely for $|\la|$ up to~$10^4$
(the coefficient~$10^{-3}$ was not determined by any computation, in
contrast to the coefficient in~\eqref{fisher-ED}, but it was chosen to
give a suitable match). When $|\la| > 10^4$, the error committed by
the Gauss--Legendre method starts to increase erratically, following
roughly the equation \mbox{$E_D = 10^{-12} \sqrt{|\la|}$}. This is
likely due to round-off error, as $|y(\xi)|$ is
approximately~$\sqrt{|\la|}$. This suggests that the loss of accuracy
in the fourth-order Magnus method is also due to round-off error,
exacerbated by ill-conditioning of the matrix exponential (compare
with the influence of the method for computing the matrix exponential,
as noted on~\pageref{expm}).

\section{Conclusions}
\label{sec:concl}

We found that the fourth-order Magnus method, when applied to the
linear differential equation~\eqref{bvp3} in the regime $|\la| \gg
1/h^2$, commits a global error of order~$h^2/\sqrt{|\la|}$ (relative
to the exact solution). It is remarkable that the Magnus method
converges at all. The convergence
result~\cite{moan.niesen:convergence} for the Magnus series mentioned
earlier guarantees convergence only when $|\la| < \pi/h^2$, so the
usual convergence proof for the truncated series does not
hold. However, the error analysis in this paper shows that the method
does indeed converge for equations of the form~\eqref{bvp3}.

Given that the method converges, it is remarkable that the order of
the method drops. This is connected to the concept of stability. The
Magnus method solves autonomous linear equations exactly.
\emph{A~fortiori}, the numerical solution of Dahlquist's test equation
$y' = ay$ (with $a\in\C$) is stable if and only if the exact solution
is stable, meaning that it converges to~0 as $x\to\infty$. Hence, the
Magnus method is A-stable and even L-stable (see, e.g., Hairer and
Wanner~\cite{hairer.wanner:solving} for a definition of these
terms). Nevertheless, the fourth-order Magnus method suffers from
order reduction in the current setting, in which the equation is
\emph{nearly} autonomous. This may be connected to the fact that the
fourth-order Magnus method is not B-stable. A simple counterexample is
given by the equation $y' = Ay$ with $A(x) = \bigl[
\begin{smallmatrix} -1 & 1 \\ 0 & -1
\end{smallmatrix} \bigr]$ for $x<3$ and $A(x) = \bigl[
\begin{smallmatrix} -1 & 0 \\ 1 & -1 
\end{smallmatrix} \bigr]$ for $x>3$. This equation is contractive, but
the numerical solution does not preserve contractivity as $h$
increases above~6. In contrast, the Gauss--Legendre method and the
exponential midpoint rule are known to be B-stable, and they do not
suffer from order reduction. We refer again to Hairer and
Wanner~\cite{hairer.wanner:solving} for a precise definition of
B-stability and its connection to order reduction.

Similar results were obtained by Hochbruck and
Lubich~\cite{hochbruck.lubich:on*1}, who treated Magnus methods
applied to semi-discretized Schr\"odinger equations. They could prove
that the method converges even when there is no known convergence
result for the untruncated Magnus series. Gonz{\'a}lez, Ostermann and
Thalhammer~\cite{gonzalez.ostermann.ea:second-order} found that the
exponential midpoint rule suffers from order reduction when applied to
semi-discretized parabolic equations. The matrices in the
semi-descretized equations considered by them have negative
eigenvalues that are large in magnitude, just as the
problem~\eqref{bvp3} treated here.

However, the fourth-order Magnus method regains the full order when
combining the solution of the differential equation~\eqref{bvp3}
satisfying the boundary condition at $\xi=-\infty$ with the one
satisfying the condition at $\xi=+\infty$ to form the Evans function.
As is clearly shown both by the analysis and by the experiment, the
error committed by the fourth-order Magnus is of order~$h^4$ uniformly
in~$\la$. Nevertheless, the Gauss--Legendre method is still superior:
its error decreases as $|\la|$ increases. 

The same holds to a lesser degree for the exponential midpoint
rule. The analysis indicates that the error commited by this method is
of order~$\la^{-1/2}h^2$. The numerical results for the exponential
midpoint rule do not quite seem to agree with this, but they also show
that the error decreases as a function of~$|\la|$; the reason for this
discrepancy is unknown. Nevertheless, we can conclude that the second
term in the Magnus expansion~\eqref{magnus} actually harms the
numerical algorithm when $\la$ is large in magnitude.

This suggests that the Right Correction Magnus Series, as proposed by
Degani and Schiff~\cite{degani.schiff:rcms}, or the modified Magnus
method, as proposed by Iserles~\cite{iserles:on*1}, might perform well
on this problem. We ran some preliminary experiments with these
methods, which showed that the error in the stiff component is greatly
reduced and comparable to the error committed by the Gauss--Legendre
method. However, the nonstiff component seems to suffer from round-off
error. A full analysis of these methods warrants further
investigation.

As explained at the start of Section~\ref{sec:evans}, our interest lies in
the stability analysis for travelling waves for the
reaction--diffusion equation~\eqref{rd}. By energy estimates similar
to those in Brin~\cite[\S3.2]{brin:numerical*1}, we find that the
eigenvalues are contained in the wedge given by
\begin{equation}
\label{eq:wedge}
\begin{gathered}
\Re\la \le \tfrac14 c^2 + \max\nolimits_\xi |f'(\hu(\xi))|, \\
\Re\la + |\Im\la | \le c^2 + \max\nolimits_\xi |f'(\hu(\xi))|. 
\end{gathered}
\end{equation}
For the Fisher equation used as an example in the previous section, we
have $\max_\xi |f'(\hu(\xi))| = 1$. Therefore, the analysis reported
in this paper is of limited use when assessing the stability of the
travelling wave for this equation. However, for other equations the
wedge~\eqref{eq:wedge} and the regime in which our analysis is valid
may overlap.

We mentioned in the introduction that this work was instigated by the
paper of Aparicio, Malham and
Oliver~\cite{aparicio.malham.ea:numerical}. That paper studies a
similar equation in the regime $\arg\la = \pi$, while the analysis in
the current paper concerns the complementary regime $|\arg \la| <
\pi-\eps$ with $\eps > 0$. In fact, it is possible to derive a
combined estimate for the local error valid in both regimes using
WKB-analysis. 

Finally, we wish to stress that the analysis reported here is only
valid for \emph{scalar} reaction--diffusion equations. We have not
done a full analysis for systems of reaction--diffusion equations of
the form~\eqref{rd}. However, numerical experiments suggest that also
in this case, the Magnus method suffers from order reduction but
recovers the full order when matching the solutions to evaluate the
Evans function. Generally, the Gauss--Legendre method still
outperforms the Magnus method when $\la$ is away from the essential
spectrum, but the difference is not so pronounced.

\subsection*{Acknowledgements}

We would like to thank Arjen Doelman, Marco Marletta, Per Christian
Moan and Khaled Saad for stimulating discussions during the
preparation of this paper. We also thank the anonymous referees for
their comments and suggestions for improvement of the manuscript.

\bibliography{../jitse}
\bibliographystyle{amsplain}

\end{document}